\documentclass[a4paper]{amsart}
\usepackage[T1]{fontenc}
\usepackage[utf8]{inputenc}
\usepackage[active]{srcltx}
\usepackage{amsthm}
\usepackage{amstext}
\usepackage{amssymb}
\usepackage{esint}
\usepackage{xargs}[2008/03/08]

\makeatletter

\newcommand{\noun}[1]{\textsc{#1}}

\theoremstyle{plain}
\newtheorem{thm}{\protect\theoremname}
  \theoremstyle{plain}
  \newtheorem{cor}[thm]{\protect\corollaryname}

\makeatother

  \providecommand{\corollaryname}{Corollary}
\providecommand{\theoremname}{Theorem}

\begin{document}

\title{Higher-Order Daehee numbers and polynomials}

\author{Dae San Kim and Taekyun Kim}
\begin{abstract}
Recently, Daehee numbers and polynomials are introduced by the authors.
In this paper, we consider the Daehee numbers and polynomials of order
$k\left(\in\mathbb{N}\right)$ and give some relation between Daehee
polynomials of order $k$$\left(\in\mathbb{N}\right)$ and special
polynomials.

\global\long\def\zp{\mathbb{Z}_{p}}

\newcommandx\ds[1][usedefault, addprefix=\global, 1=]{\widehat{D}_{n}^{#1}}

\end{abstract}
\maketitle

\section{Introduction}

For $\alpha\in\mathbb{N}$, as is well known, the Bernoulli
polynomials of order $\alpha$ are defined by the generating function to be
\begin{equation}
\left(\frac{t}{e^{t}-1}\right)^{\alpha}e^{xt}=\sum_{n=0}^{\infty}B_{n}^{\left(\alpha\right)}\left(x\right)\frac{t^{n}}{n!},\label{eq:1}
\end{equation}
(see {[}1-14{]}).

When $x=0$, $B_{n}^{\left(\alpha\right)}=B_{n}^{\left(\alpha\right)}\left(0\right)$
are the Bernoulli numbers of order $\alpha$. In \cite{2013DSKIM,MR2390695,MR2479746},
the Daehee polynomials are defined by the generating function
to be
\begin{equation}
\left(\frac{\log\left(1+t\right)}{t}\right)\left(1+t\right)^{x}=\sum_{n=0}^{\infty}D_{n}\left(x\right)\frac{t^{n}}{n!}.\label{eq:2}
\end{equation}

When $x=0$, $D_{n}=D_{n}\left(0\right)$ are called the Daehee
numbers.

Throughout this paper, $\mathbb{Z}_{p}$, $\mathbb{Q}_{p}$ and $\mathbb{C}_{p}$
will denote the ring of $p$-adic integers, the field of $p$-adic
numbers and the completion of algebraic closure of $\mathbb{Q}_{p}$.
The $p$-adic norm $\left|\cdot\right|_{p}$ is normalized as $\left|p\right|_{p}=\frac{1}{p}$.
Let $\textnormal{UD}\left(\mathbb{Z}_{p}\right)$ be the space of
uniformly differentiable functions on $\mathbb{Z}_{p}$. For $f\in\textnormal{UD}\left(\mathbb{Z}_{p}\right)$,
the $p$-adic invariant integral on $\mathbb{Z}_{p}$ is defined by
\begin{equation}
I\left(f\right)=\int_{\mathbb{Z}_{p}}f\left(x\right)d\mu\left(x\right)=\lim_{n\rightarrow\infty}\frac{1}{p^{n}}\sum_{x=0}^{p^{n}-1}f\left(x\right),\label{eq:3}
\end{equation}
(see \cite{MR2845943}).

Let $f_{1}\left(x\right)=f\left(x+1\right).$ Then, by (\ref{eq:3}),
we get
\begin{equation}
I\left(f_{1}\right)-I\left(f\right)=f^{\prime}\left(0\right),\textrm{ where }f^{\prime}\left(0\right)=\left.\frac{df\left(x\right)}{dx}\right|_{x=0}.\label{eq:4}
\end{equation}

The signed Stirling numbers of the first kind $S_{1}(n,l)$ are defined by
\begin{eqnarray}
\left(x\right)_{n} & = & x\left(x-1\right)\cdots\left(x-n+1\right)\label{eq:5}\\
 & = & \sum_{l=0}^{\infty}S_{1}\left(n,l\right)x^{l},\nonumber
\end{eqnarray}
(see \cite{MR2508979,MR2597988,2013DSKIM}).

From (\ref{eq:5}), we note that
\begin{eqnarray*}
x^{\left(n\right)} & = & x\left(x+1\right)\cdots\left(x+n-1\right)=\left(-1\right)^{n}\left(-x\right)_{n}\\
 & = & \sum_{l=0}^{n}\left(-1\right)^{n-l}S_{1}\left(n,l\right)x^{l},
\end{eqnarray*}
(see \cite{2013DSKIM,MR2931605,MR2479746}).

The Stirling numbers of the second kind $S_{2}(l,n)$ are defined by the generating
function to be
\begin{eqnarray}
\left(e^{t}-1\right)^{n} & = & n!\sum_{l=n}^{\infty}S_{2}\left(l,n\right)\frac{t^{l}}{l!}\label{eq:6}\\
 & = & \sum_{l=0}^{\infty}\frac{n!}{\left(l+n\right)!}S_{2}\left(l+n,n\right)t^{l+n}.\nonumber
\end{eqnarray}

In this paper, we study the higher-order Daehee numbers and polynomials
and give some relations between Daehee polynomials and special polynomials.

\section{Higher-order Daehee polynomials}

In this section, we assume that $t\in\mathbb{C}_{p}$ with $\left|t\right|_{p}<p^{\frac{-1}{p-1}}$.

For $k\in\mathbb{N}$, let us consider the Daehee numbers of the first
kind of order $k$ :

\begin{equation}
D_{n}^{\left(k\right)}=\underset{k-\textrm{times}}{\underbrace{\int_{\mathbb{Z}_{p}}\cdots\int_{\mathbb{Z}_{p}}}}\left(x_{1}+x_{2}+\cdots+x_{k}\right)_{n}d\mu\left(x_{1}\right)\cdots d\mu\left(x_{k}\right),\label{eq:7}
\end{equation}
where $n\in\mathbb{Z}_{\ge0}$.

From (\ref{eq:7}), we can derive the generating function of $D_{n}^{\left(k\right)}$
as follows :
\begin{align}
 & \sum_{n=0}^{\infty}D_{n}^{\left(k\right)}\frac{t^{n}}{n!}\label{eq:8}\\
= & \int_{\zp}\cdots\int_{\zp}\sum_{n=0}^{\infty}\dbinom{x_{1}+\cdots+x_{k}}{n}t^{n}d\mu\left(x_{1}\right)\cdots d\mu\left(x_{k}\right)\nonumber \\
= & \int_{\mathbb{Z}_{p}}\cdots\int_{\mathbb{Z}_{p}}\left(1+t\right)^{x_{1}+\cdots+x_{k}}d\mu\left(x_{1}\right)\cdots d\mu\left(x_{k}\right).\nonumber
\end{align}

By (\ref{eq:4}), we easily see that
\begin{equation}
\int_{\zp}\left(1+t\right)^{x}d\mu\left(x\right)=\frac{\log\left(1+t\right)}{t}.\label{eq:9}
\end{equation}

Thus, by (\ref{eq:8}) and (\ref{eq:9}), we get
\begin{equation}
\sum_{n=0}^{\infty}D_{n}^{\left(k\right)}\frac{t^{n}}{n!}=\left(\frac{\log\left(1+t\right)}{t}\right)^{k}.\label{eq:10}
\end{equation}

Now, we observe that
\begin{eqnarray}
\left(\frac{\log\left(1+t\right)}{t}\right)^{k} & = & \frac{k!}{t^{k}}\sum_{l=k}^{\infty}S_{1}\left(t,k\right)\frac{t^{l}}{l!}\label{eq:11}\\
 & = & \sum_{n=0}^{\infty}S_{1}\left(n+k,k\right)\frac{k!}{\left(n+k\right)!}t^{n}\nonumber \\
 & = & \sum_{n=0}^{\infty}\frac{S_{1}\left(n+k,k\right)}{\tbinom{n+k}{k}}\frac{t^{n}}{n!}.\nonumber
\end{eqnarray}

Therefore, by (\ref{eq:10}) and (\ref{eq:11}), we obtain the following
theorem.
\begin{thm}
For $n\in\mathbb{Z}_{\ge0}$, $k\in\mathbb{N}$, we have
\[
D_{n}^{\left(k\right)}=\frac{S_{1}\left(n+k,k\right)}{\tbinom{n+k}{k}}.
\]

\end{thm}

It is easy to show that
\begin{equation}
\left(\frac{\log\left(1+t\right)}{t}\right)^{k}=\sum_{n=0}^{\infty}B_{n}^{\left(n+k+1\right)}\left(1\right)\frac{t^{n}}{n!}.\label{eq:12}
\end{equation}

Threfore, we obtain the following corollary.
\begin{cor}
For $n\in\mathbb{Z}_{\ge0}$, $k\in\mathbb{N}$, we have
\[
D_{n}^{\left(k\right)}=\frac{S_{1}\left(n+k,k\right)}{\tbinom{n+k}{k}}=B_{n}^{\left(n+k+1\right)}\left(1\right).
\]

\end{cor}
From (\ref{eq:7}), we note that
\begin{eqnarray}
D_{n}^{\left(k\right)} & = & \int_{\mathbb{Z}_{p}}\cdots\int_{\mathbb{Z}_{p}}\left(x_{1}+\cdots+x_{k}\right)_{n}d\mu\left(x_{1}\right)\cdots d\mu\left(x_{k}\right)\label{eq:13}\\
 & = & \sum_{l=0}^{n}S_{1}\left(n,l\right)\int_{\mathbb{Z}_{p}}\cdots\int_{\mathbb{Z}_{p}}\left(x_{1}+\cdots+x_{k}\right)^{l}d\mu\left(x_{1}\right)\cdots d\mu\left(x_{k}\right)\nonumber \\
 & = & \sum_{l=0}^{n}S_{1}\left(n,l\right)B_{l}^{\left(k\right)}.\nonumber
\end{eqnarray}

Therefore, by (\ref{eq:13}), we obtain the following theorem.
\begin{thm}
For $n\in\mathbb{Z}_{\ge0}$, $k\in\mathbb{N}$, we have
\begin{eqnarray*}
D_{n}^{\left(k\right)} & = & \sum_{l_{1}+\cdots+l_{k}=n}\dbinom{n}{l_{1},\cdots,l_{k}}D_{l_{1}}\cdots D_{l_{k}}\\
 & = & \sum_{l=0}^{n}S_{1}\left(n,l\right)B_{l}^{\left(k\right)}.
\end{eqnarray*}

\end{thm}

From (\ref{eq:10}), we can derive
\begin{equation}
\sum_{n=0}^{\infty}D_{n}^{\left(k\right)}\frac{\left(e^{t}-1\right)^{n}}{n!}=\left(\frac{t}{e^{t}-1}\right)^{k}=\sum_{n=0}^{\infty}B_{n}^{\left(k\right)}\frac{t^{n}}{n!},\label{eq:14}
\end{equation}
 and
\begin{equation}
\sum_{n=0}^{\infty}D_{n}^{\left(k\right)}\frac{\left(e^{t}-1\right)^{n}}{n!}=\sum_{m=0}^{\infty}\left(\sum_{n=0}^{m}D_{n}^{\left(k\right)}S_{2}\left(n,m\right)\right)\frac{t^{n}}{m!}.\label{eq:15}
\end{equation}

Therefore, by (\ref{eq:14}) and (\ref{eq:15}), we obtain the following
theorem.
\begin{thm}
For $m\in\mathbb{Z}_{\ge0}$, $k\in\mathbb{N}$, we have
\[
B_{m}^{\left(k\right)}=\sum_{n=0}^{m}D_{n}^{\left(k\right)}S_{2}\left(m,n\right).
\]

\end{thm}

Now, we consider the higher-order Daehee polynomials as follows :

\begin{eqnarray}
D_{n}^{\left(k\right)}\left(x\right) & = & \int_{\mathbb{Z}_{p}}\cdots\int_{\mathbb{Z}_{p}}\left(x_{1}+\cdots+x_{k}+x\right)_{n}d\mu\left(x_{1}\right)\cdots d\mu\left(x_{k}\right)\label{eq:16}.
\end{eqnarray}

Thus, by (\ref{eq:16}), we get
\begin{align}
 & D_{n}^{\left(k\right)}\left(x\right)\label{eq:17}\\
= & \sum_{l=0}^{n}S_{1}\left(n,l\right)\int_{\mathbb{Z}_{p}}\cdots\int_{\mathbb{Z}_{p}}\left(x_{1}+\cdots+x_{k}+x\right)^{l}d\mu\left(x_{1}\right)\cdots d\mu\left(x_{k}\right)\nonumber \\
= & \sum_{l=0}^{n}S_{1}\left(n,l\right)B_{l}^{\left(k\right)}\left(x\right).\nonumber
\end{align}

Therefore, by (\ref{eq:17}), we obtain the following theorem.
\begin{thm}
For $n\in\mathbb{Z}_{\ge0}$, $k\in\mathbb{N}$, we have
\[
D_{n}^{\left(k\right)}\left(x\right)=\sum_{l=0}^{n}S_{1}\left(n,l\right)B_{l}^{\left(k\right)}\left(x\right).
\]

\end{thm}

From (\ref{eq:16}), we derive the generating function of $D_{n}^{\left(k\right)}\left(x\right)$:

\begin{align}
 & \sum_{n=0}^{\infty}D_{n}^{(k)}\left(x\right)\frac{t^{n}}{n!}\label{eq:18}\\
= & \int_{\mathbb{Z}_{p}}\cdots\int_{\mathbb{Z}_{p}}\sum_{n=0}^{\infty}\dbinom{x_{1}+\cdots+x_{k}+x}{n}t^{n}d\mu\left(x_{1}\right)\cdots d\mu\left(x_{k}\right)\nonumber \\
= & \int_{\mathbb{Z}_{p}}\cdots\int_{\mathbb{Z}_{p}}\left(1+t\right)^{x_{1}+\cdots+x_{k}+x}d\mu\left(x_{1}\right)\cdots d\mu\left(x_{k}\right)\nonumber \\
= & \left(\frac{\log\left(1+t\right)}{t}\right)^{k}\left(1+t\right)^{x}.\nonumber
\end{align}

It is easy to show that
\begin{equation}
\left(\frac{\log\left(1+t\right)}{t}\right)^{k}\left(1+t\right)^{x}=\sum_{n=0}^{\infty}B_{n}^{\left(n+k+1\right)}\left(x+1\right)\frac{t^{n}}{n!}.\label{eq:19}
\end{equation}

Therefore, by (\ref{eq:18}) and (\ref{eq:19}), we obtain the following
theorem.
\begin{thm}
For $n\in\mathbb{Z}_{\ge0}$, $k\in\mathbb{N}$,
\begin{eqnarray*}
D_{n}^{\left(k\right)}\left(x\right) & = & B_{n}^{\left(n+k+1\right)}\left(x+1\right)\\
 & = & \sum_{l=0}^{n}\dbinom{n}{l}B_{l}^{\left(n+k+1\right)}\left(x+1\right)^{n-l}.
\end{eqnarray*}

\end{thm}

In (\ref{eq:18}), we note that
\begin{equation}
\sum_{n=0}^{\infty}D_{n}^{\left(k\right)}\left(x\right)\frac{\left(e^{t}-1\right)^{n}}{n!}=\sum_{m=0}^{\infty}\left(\sum_{n=0}^{m}S_{2}\left(n,m\right)D_{n}^{\left(k\right)}\left(x\right)\right)\frac{t^{m}}{m!}\label{eq:20}
\end{equation}
and
\begin{eqnarray}
\sum_{n=0}^{\infty}D_{n}^{\left(k\right)}\left(x\right)\frac{\left(e^{t}-1\right)^{n}}{n!} & = & \left(\frac{t}{e^{t}-1}\right)^{k}e^{xt}\label{eq:21}\\
 & = & \sum_{m=0}^{\infty}B_{m}^{\left(k\right)}\left(x\right)\frac{t^{m}}{m!}.\nonumber
\end{eqnarray}

Therefore, by (\ref{eq:20}) and (\ref{eq:21}), we obtain the following
theorem.
\begin{thm}
For $m\in\mathbb{Z}_{\ge0}$, $k\in\mathbb{N}$, we have
\[
B_{m}^{\left(k\right)}\left(x\right)=\sum_{n=0}^{m}S_{2}\left(m,n\right)D_{n}^{\left(k\right)}\left(x\right).
\]

\end{thm}

Now, we define Daehee numbers of the second kind of order $k$$\left(\in\mathbb{N}\right)$
:
\begin{align}
 & \ds[\left(k\right)]\label{eq:22}\\
= & \left(-1\right)^{n}\int_{\zp}\cdots\int_{\zp}\left(-x_{1}-x_{2}-\cdots-x_{k}\right)_{n}d\mu\left(x_{1}\right)\cdots d\mu\left(x_{k}\right)\nonumber \\
= & \left(-1\right)^{n}\sum_{l=0}^{n}\left(-1\right)^{n-l}S_{1}\left(n,l\right)B_{l}^{\left(k\right)}=\sum_{l=0}^{n}\begin{bmatrix}n\\
l
\end{bmatrix}B_{l}^{\left(k\right)},\nonumber
\end{align}
where $\begin{bmatrix}n\\
l
\end{bmatrix}=\left(-1\right)^{n-l}S_{1}$$\left(n,l\right)$.

Thus, by (\ref{eq:22}), we get
\begin{align}
 & \ds[\left(k\right)]\label{eq:23}\\
= & \left(-1\right)^{n}\int_{\zp}\cdots\int_{\zp}\left(-x_{1}-x_{2}-\cdots-x_{k}\right)_{n}d\mu\left(x_{1}\right)\cdots d\mu\left(x_{k}\right)\nonumber \\
= & \left(-1\right)^{n}\sum_{l=0}^{n}S_{1}\left(n,l\right)\left(-1\right)^{l}\int_{\zp}\cdots\int_{\zp}\left(x_{1}+x_{2}+\cdots+x_{k}\right)^{l}d\mu\left(x_{1}\right)\cdots d\mu\left(x_{k}\right)\nonumber \\
= & \sum_{l=0}^{n}\left(-1\right)^{n-l}S_{1}\left(n,l\right)B_{l}^{\left(k\right)}=\sum_{l=0}^{n}\begin{bmatrix}n\\
l
\end{bmatrix}B_{l}^{\left(k\right)},\nonumber
\end{align}
where $\begin{bmatrix}n\\
l
\end{bmatrix}=\left(-1\right)^{n-l}S_{1}$$\left(n,l\right)$.

Therefore, by (\ref{eq:23}), we obtain the following theorem.
\begin{thm}
For $n\in\mathbb{Z}_{\ge0}$, $k\in\mathbb{N}$, we have
\[
\ds[\left(k\right)]=\sum_{l=0}^{n}\begin{bmatrix}n\\
l
\end{bmatrix}B_{l}^{\left(k\right)}.
\]

\end{thm}

From (\ref{eq:22}), we derive the generating function of $\ds[\left(k\right)]$:

\begin{align}
 & \sum_{n=0}^{\infty}\ds[\left(k\right)]\frac{t^{n}}{n!}\label{eq:24}\\
= & \int_{\zp}\cdots\int_{\zp}\sum_{n=0}^{\infty}\dbinom{x_{1}+\cdots+x_{k}+n-1}{n}t^{n}d\mu\left(x_{1}\right)\cdots d\mu\left(x_{k}\right)\nonumber \\
= & \int_{\mathbb{Z}_{p}}\cdots\int_{\mathbb{Z}_{p}}\left(1-t\right)^{-x_{1}-\cdots-x_{k}}d\mu\left(x_{1}\right)\cdots d\mu\left(x_{k}\right)\nonumber \\
= & \left(\frac{\left(1-t\right)\log\left(1-t\right)}{-t}\right)^{k}.\nonumber
\end{align}

By (\ref{eq:24}), we get
\begin{eqnarray}
\sum_{n=0}^{\infty}\ds[\left(k\right)]\frac{\left(1-e^{-t}\right)^{n}}{n!} & = & \left(\frac{e^{-t}\left(-t\right)}{e^{-t}-1}\right)^{k}=\left(\frac{t}{e^{t}-1}\right)^{k}\label{eq:25}\\
 & = & \sum_{m=0}^{\infty}B_{m}^{\left(k\right)}\frac{t^{m}}{m!},\nonumber
\end{eqnarray}
and
\begin{equation}
\sum_{n=0}^{\infty}\ds[\left(k\right)]\frac{\left(1-e^{-t}\right)^{n}}{n!}=\sum_{m=0}^{\infty}\left(\sum_{n=0}^{m}\ds[\left(k\right)]\left(-1\right)^{m-n}S_{2}\left(m,n\right)\right)\frac{t^{m}}{m!}.\label{eq:26}
\end{equation}

Thererfore, by (\ref{eq:25}) and (\ref{eq:26}), we obtain the following
theorem.
\begin{thm}
For $m\in\mathbb{Z}_{\ge0}$, $k\in\mathbb{N}$, we have
\[
B_{m}^{\left(k\right)}=\sum_{n=0}^{m}\ds[\left(k\right)]\left(-1\right)^{n-m}S_{2}\left(m,n\right).
\]

\end{thm}

Now, we consider the higher-order Daehee polynomials of the second
kind :
\begin{equation}
\ds[\left(k\right)]\left(x\right)=\int_{\zp}\cdots\int_{\zp}\left(x_{1}+x_{2}+\cdots+x_{k}-x\right)^{(n)}d\mu\left(x_{1}\right)\cdots d\mu\left(x_{k}\right).\label{eq:27}
\end{equation}

Thus, by (\ref{eq:27}), we get
\begin{align}
 & \ds[\left(k\right)]\left(x\right)\label{eq:28}\\
= & \left(-1\right)^{n}\int_{\zp}\cdots\int_{\zp}\left(-x_{1}-x_{2}-\cdots-x_{k}+x\right)_{n}d\mu\left(x_{1}\right)\cdots d\mu\left(x_{k}\right)\nonumber \\
= & \left(-1\right)^{n}\sum_{l=0}^{n}S_{1}\left(n,l\right)\int_{\zp}\cdots\int_{\zp}\left(-x_{1}-x_{2}-\cdots-x_{k}+x\right)^{l}d\mu\left(x_{1}\right)\cdots d\mu\left(x_{k}\right)\nonumber \\
= & \left(-1\right)^{n}\sum_{l=0}^{n}S_{1}\left(n,l\right)\sum_{m=0}^{l}\dbinom{l}{m}x^{l-m}\int_{\zp}\cdots\int_{\zp}\left(-x_{1}-x_{2}-\cdots-x_{k}\right)^{m}d\mu\left(x_{1}\right)\cdots d\mu\left(x_{k}\right)\nonumber \\
= & \left(-1\right)^{n}\sum_{l=0}^{n}S_{1}\left(n,l\right)\sum_{m=0}^{l}\dbinom{l}{m}\left(-1\right)^{m}x^{l-m}B_{m}^{\left(k\right)}\nonumber \\
= & \sum_{l=0}^{n}\left(-1\right)^{n-l}S_{1}\left(n,l\right)B_{l}^{\left(k\right)}\left(-x\right).\nonumber
\end{align}

Thus, by (\ref{eq:28}), we get
\begin{equation}
\ds[\left(k\right)]\left(x\right)=\sum_{l=0}^{n}\left(-1\right)^{n-l}S_{1}\left(n,l\right)B_{l}^{\left(k\right)}\left(-x\right).\label{eq:29}
\end{equation}

Let us consider the generating function of $D_{n}^{\left(k\right)}\left(x\right)$
as follows :
\begin{align}
 & \sum_{n=0}^{\infty}\ds[\left(k\right)]\left(x\right)\frac{t^{n}}{n!}\label{eq:30}\\
= & \int_{\zp}\cdots\int_{\zp}\sum_{n=0}^{\infty}\dbinom{x_{1}+\cdots+x_{k}-x+n-1}{n}t^{n}d\mu\left(x_{1}\right)\cdots d\mu\left(x_{k}\right)\nonumber \\
= & \int_{\zp}\cdots\int_{\zp}\left(1-t\right)^{-x_{1}-\cdots-x_{k}+x}d\mu\left(x_{1}\right)\cdots d\mu\left(x_{k}\right)\nonumber \\
= & \left(\frac{\left(1-t\right)\log\left(1-t\right)}{-t}\right)^{k}\left(1-t\right)^{x}.\nonumber
\end{align}

From (\ref{eq:30}), we have

\begin{align}
 & \sum_{n=0}^{\infty}\ds[\left(k\right)]\left(x\right)\left(-1\right)^{n}\frac{t^{n}}{n!}\label{eq:31}\\
= & \left(\frac{\log\left(1+t\right)}{t}\right)^{k}\left(1+t\right)^{x+k}\nonumber \\
= & \sum_{n=0}^{\infty}B_{n}^{\left(n+k+1\right)}\left(x+k+1\right)\frac{t^{n}}{n!}.\nonumber
\end{align}

Therefore, by (\ref{eq:31}), we obtain the following theorem.
\begin{thm}
For $n\in\mathbb{Z}_{\ge0}$, $k\in\mathbb{N}$, we have
\[
\left(-1\right)^{n}\ds[\left(k\right)]\left(x\right)=B_{n}^{\left(n+k+1\right)}\left(x+k+1\right).
\]

\end{thm}

By (\ref{eq:30}), we get
\begin{eqnarray}
\sum_{n=0}^{\infty}\ds[\left(k\right)]\left(x\right)\frac{\left(1-e^{-t}\right)^{n}}{n!} & = & e^{-tx}\left(\frac{t}{e^{t}-1}\right)^{k}\label{eq:32}\\
 & = & \sum_{m=0}^{\infty}B_{m}^{\left(k\right)}\left(-x\right)\frac{t^{m}}{m!},\nonumber
\end{eqnarray}
and
\begin{align}
 & \sum_{n=0}^{\infty}\ds[\left(k\right)]\left(x\right)\frac{1}{n!}\left(1-e^{-t}\right)^{n}\label{eq:33}\\
= & \sum_{m=0}^{\infty}\left(\sum_{n=0}^{m}\ds[\left(k\right)]\left(x\right)\left(-1\right)^{m-n}S_{2}\left(m,n\right)\right)\frac{t^{m}}{m!}.\nonumber
\end{align}

Therefore, by (\ref{eq:32}) and (\ref{thm:12}), we obtain the following
theorem.
\begin{thm}
For $m\in\mathbb{Z}_{\ge0}$, $k\in\mathbb{N}$, we have
\[
B_{m}^{\left(k\right)}\left(-x\right)=\sum_{n=0}^{m}\ds[\left(k\right)]\left(x\right)\left(-1\right)^{m-n}S_{2}\left(m,n\right).
\]

\end{thm}

Now, we observe that
\begin{align}
 & \left(-1\right)^{n}\frac{D_{n}^{\left(k\right)}\left(x\right)}{n!}\label{eq:34}\\
= & \left(-1\right)^{n}\int_{\zp}\cdots\int_{\zp}\dbinom{x_{1}+\cdots+x_{k}+x}{n}d\mu\left(x_{1}\right)\cdots d\mu\left(x_{k}\right)\nonumber \\
= & \int_{\zp}\cdots\int_{\zp}\dbinom{-(x_{1}+\cdots+x_{k})-x+n-1}{n}d\mu\left(x_{1}\right)\cdots d\mu\left(x_{k}\right)\nonumber \\
= & \sum_{m=0}^{n}\dbinom{n-1}{n-m}\int_{\zp}\cdots\int_{\zp}\dbinom{-(x_{1}+\cdots+x_{k})-x}{m}d\mu\left(x_{1}\right)\cdots d\mu\left(x_{k}\right)\nonumber \\
= & \sum_{m=0}^{n}\frac{\tbinom{n-1}{n-m}}{m!}m!\int_{\zp}\cdots\int_{\zp}\dbinom{-(x_{1}+\cdots+x_{k})-x}{m}d\mu\left(x_{1}\right)\cdots d\mu\left(x_{k}\right)\nonumber \\
= & \sum_{m=1}^{n}\frac{\tbinom{n-1}{n-m}}{m!}\left(-1\right)^{m}\widehat{D}_{m}^{\left(k\right)}\left(-x\right).\nonumber
\end{align}

Therefore, by (\ref{eq:34}), we obtain the following theorem.
\begin{thm}
\label{thm:12}For $n\in\mathbb{Z}_{\ge0}$, $k\in\mathbb{N}$, we
have
\[
\left(-1\right)^{n}\frac{D_{n}^{\left(k\right)}\left(x\right)}{n!}=\sum_{m=1}^{n}\frac{\tbinom{n-1}{n-m}}{m!}\left(-1\right)^{m}\widehat{D}_{m}^{\left(k\right)}\left(-x\right).
\]

\end{thm}

By the same method as Theorem \ref{thm:12}, we get
\begin{align}
 & \frac{\ds[\left(k\right)]\left(x\right)}{n!}\label{eq:35}\\
= & \int_{\zp}\cdots\int_{\zp}\dbinom{x_{1}+\cdots+x_{k}-x+n-1}{n}d\mu\left(x_{1}\right)\cdots d\mu\left(x_{k}\right)\nonumber \\
= & \sum_{m=0}^{n}\dbinom{n-1}{n-m}\int_{\zp}\cdots\int_{\zp}\dbinom{x_{1}+\cdots+x_{k}-x}{m}d\mu\left(x_{1}\right)\cdots d\mu\left(x_{k}\right)\nonumber \\
= & \sum_{m=0}^{n}\frac{\tbinom{n-1}{n-m}}{m!}m!\int_{\zp}\cdots\int_{\zp}\dbinom{x_{1}+\cdots+x_{k}-x}{m}d\mu\left(x_{1}\right)\cdots d\mu\left(x_{k}\right)\nonumber \\
= & \sum_{m=1}^{n}\frac{\tbinom{n-1}{n-m}}{m!}D_{m}^{\left(k\right)}\left(-x\right).\nonumber
\end{align}

Thus, by (\ref{eq:35}), we get
\begin{equation}
\frac{\ds[\left(k\right)]\left(x\right)}{n!}=\sum_{m=1}^{n}\frac{\tbinom{n-1}{n-m}}{m!}D_{m}^{\left(k\right)}\left(-x\right).\label{eq:36}
\end{equation}

\bibliographystyle{plain}
\nocite{*}
\bibliography{higher-daehee}

$\,$

\noindent \noun{Department of Mathematics, Sogang University, Seoul
121-742, Republic of Korea}

\noindent \emph{E-mail}\noun{ }\emph{address : }\texttt{dskim@sogang.ac.kr}\\

\noindent \noun{Department of Mathematics, Kwangwoon University, Seoul
139-701, Republic of Korea}

\noindent \emph{E-mail}\noun{ }\emph{address : }\texttt{tkkim@kw.ac.kr}
\end{document}